\documentclass[12pt]{article}
\usepackage{epsfig}
\usepackage{amsmath}
\usepackage{amsthm}
\usepackage{amsfonts}
\usepackage[english]{babel}
\usepackage{amssymb}
\usepackage[T1]{fontenc}
\usepackage{amsmath}

\newtheorem{Th}{Theorem}
\newtheorem{Def}{Definition}
\newtheorem{lemma}{Lemma}
\newtheorem{example}{Example}

\newtheorem{Prop}{Proposition}

\newcommand{\A}{\mathcal{A}}

\newcommand{\E}{\mathbb{E}}
\newcommand{\F}{\mathcal{F}}
\newcommand{\R}{\mathbb{R}}
\newcommand{\Z}{\mathbb{Z}}
\newcommand{\N}{\mathbb{N}}

\newcommand{\Q}{\mathcal{Q}}
\renewcommand{\P}{\mathbb{P}}
\newcommand{\ds}{\displaystyle}

\newcounter{tictac}

\def\1{\,\rlap{\mbox{\small\rm 1}}\kern.15em 1}
\def\ind#1{\1_{#1}}
\def\build#1_#2^#3{\mathrel{\mathop{\kern 0pt#1}\limits_{#2}^{#3}}}
\def\tend#1#2{\build\hbox to 12mm{\rightarrowfill}_{#1\rightarrow #2}^{a.s.}}

\def\converge#1#2#3{\build\hbox to
15mm{\rightarrowfill}_{#1\rightarrow #2}^{\hbox{\scriptsize #3}}}
\begin{document}

\begin{center}
{\large
    {\sc
A central limit theorem for stationary random fields
    }
}
\bigskip

Mohamed EL MACHKOURI, Dalibor VOLN\'Y

\medskip
{\it

Laboratoire de Math\'ematiques Rapha\"el Salem\\
UMR CNRS 6085, Universit\'e de Rouen (France)\\
}
\medskip
and\\
\medskip
Wei Biao WU

\medskip
{\it
University of Chicago
}
\end{center}
\centerline{\today}

{\renewcommand\abstractname{Abstract}
\begin{abstract}
\baselineskip=18pt This paper establishes a central limit theorem
and an invariance principle for a wide class of stationary random
fields under natural and easily verifiable conditions. More
precisely, we deal with random fields of the form $X_k =
g\left(\varepsilon_{k-s}, s \in \Z^d \right)$, $k\in\Z^d$, where
$(\varepsilon_i)_{i\in\Z^d}$ are i.i.d random variables and $g$ is
a measurable function. Such kind of spatial processes provides a
general framework for stationary ergodic random fields. Under a
short-range dependence condition, we show that the central limit
theorem holds without any assumption on the underlying domain on
which the process is observed. A limit theorem for the sample
auto-covariance function is also established.
\end{abstract}
\thispagestyle{empty}
\baselineskip=18pt

\section{Introduction}
Central limit theory plays a fundamental role in statistical
inference of random fields. There have been a substantial
literature for central limit theorems of random fields under
various dependence conditions. See \cite{Bolthausen1982a}, \cite{Bradley1992}, \cite{Bulinski2010}, \cite{Bulinskii1988}, \cite{Chen},
\cite{Dedecker1998},  \cite{Guyon--Richardson1984}, \cite{Jenish--Prucha2009}, \cite{Maltz1999}, \cite{Nahapetian1987}, 
\cite{Neaderhouser1978a}, \cite{Neaderhouser1978b}, \cite{Paulauskas2010}, \cite{Perera1997}, among others. However, many of them require that the
underlying random fields have very special structures such as
Gaussian, linear, Markovian or strong mixing of various types. In
applications those structural assumptions can be violated, or not
easily verifiable. 

In this paper we consider stationary random fields which are viewed as nonlinear
transforms of independent and identically distributed (iid) random
variables. Based on that representation we introduce dependence
measures and establish a central limit theorem and an invariance
principle. We assume that the random field $(X_i)_{ i \in \Z^d}$
has the form
\begin{equation}\label{definition_champ}
X_i = g\left(\varepsilon_{i-s};\,s\in\Z^d\right),
 \quad i\in\Z^d,
\end{equation}
where $(\varepsilon_j)_{j\in\Z^d}$ are iid random variables and
$g$ is a measurable function. In the one-dimensional case ($d=1$) the model 
(\ref{definition_champ}) is well known and includes linear as well as many widely
used nonlinear time series models as special cases. In Section
\ref{eq:ex} based on (\ref{definition_champ}) we shall introduce
dependence measures. It turns out that, with our dependence
measure, central limit theorems and moment inequalities can be
established in a very elegant and natural way. 

The rest of the paper is organized as follows. In Section
\ref{sec:main} we present a central limit theorem and an
invariance principle for
\begin{eqnarray*}
S_\Gamma = \sum_{i\in\Gamma}X_i,
\end{eqnarray*}
where $\Gamma$ is a finite subset of $\Z^d$ which grows to
infinity. The proof of our Theorem \ref{tlc} is based on a central limit theorem 
for $m_n$-dependent random fields established by Heinrich \cite{Heinrich88}. 
Unlike most existing results on central limit theorems
for random fields which require certain regularity conditions on
the boundary of $\Gamma$, Heinrich's central limit theorem (and consequently our Theorem \ref{tlc}) has the very
interesting property that no condition on the boundary of $\Gamma$
is needed, and the central limit theorem holds under the minimal
condition that $|\Gamma| \to \infty$, where $|\Gamma|$ the
cardinal of $\Gamma$. This is a very attractive property in
spatial applications in which the underlying observation domains
can be quite irregular. As an application, we establish a central
limit theorem for sample auto-covariances. Section \ref{sec:main}
also present an invariance principle. Proofs are provided in
Section \ref{sec:proof}.

\section{Examples and Dependence Measures}
\label{eq:ex} In (\ref{definition_champ}), we can interpret
$(\varepsilon_s)_{s \in \Z^d}$ as the input random field, $g$ is a
transform or map and $(X_i)_{i \in \Z^d}$ as the output random
field. Based on this interpretation, we define dependence measure
as follows: let $(\varepsilon_j^{'})_{j\in\Z^d}$ be an iid copy of
$(\varepsilon_j)_{j\in\Z^d}$ and consider for any positive integer
$n$ the coupled version $X_i^{\ast}$ of $X_i$ defined by
$$
X_i^{\ast}=g\left(\varepsilon^{\ast}_{i-s}\,;\,s\in\Z^d\right),
$$
where for any $j$ in $\Z^d$,
$$
\varepsilon_j^{\ast}=   \left\{\begin{array}{ll}
                        \varepsilon_j & \textrm{if $j\neq 0$} \\
            \varepsilon_0^{'}    & \textrm{if $j=0$.}
                        \end{array}\right.
$$
Recall that a Young function $\psi$ is a real convex nondecreasing
function defined on $\R^{+}$ which satisfies
$\lim_{t\to\infty}\psi(t)=\infty$ and $\psi(0)=0$.
We define the Orlicz space $\mathbb{L}_{\psi}$ as the space of real random
variables $Z$ defined on the probability space $(\Omega, \F, \P)$
such that $\E[\psi(\vert Z\vert/c)]<+\infty$ for some $c>0$. The
Orlicz space $\mathbb{L}_{\psi}$ equipped with the so-called Luxemburg
norm $\| . \|_{\psi}$ defined for any real random variable $Z$ by
\begin{eqnarray*}
\| Z\|_{\psi}=\inf\{\,c>0\,;\,\E[\psi(\vert Z\vert/c)]\leq 1\,\}
\end{eqnarray*}
is a Banach space. For more about Young functions and Orlicz
spaces one can refer to Krasnosel'skii and Rutickii \cite{K-R}.

Following Wu \cite{Wu2005}, we introduce the following dependence measures which are directly related to the underlying processes.

\begin{Def}[Physical dependence measure]
Let $\psi$ be a Young function and $i$ in $\Z^d$ be fixed. If
$X_i$ belongs to $\mathbb{L}_{\psi}$, we define the physical
dependence measure $\delta_{i,\psi}$ by
$$
\delta_{i,\psi}=\|X_i-X_i^{\ast}\|_{\psi}.
$$
If $p\in ]0,+\infty]$ and $X_i$ belongs to  $\mathbb{L}^p$, we
denote $\delta_{i,p}=\|X_i-X_i^{\ast}\|_p$.
\end{Def}
\begin{Def}[Stability]
We say that the random field $X$ defined by
$(\ref{definition_champ})$ is $p$-stable if
$$
\Delta_p:=\sum_{i\in\Z^d}\delta_{i,p}<\infty.
$$
\end{Def}
As an illustration, we give some examples of $p$-stable spatial
processes.

\begin{example}
{\em (Linear random fields) Let $(\varepsilon_i)_{i\in\Z^d}$ be
i.i.d random variables with $\varepsilon_i$ in $\mathbb{L}^p$, $p
\geq 2$. The linear random field $X$ defined for any $k$ in $\Z^d$
by
$$
X_k=\sum_{s\in\Z^d}a_s\varepsilon_{k-s}
$$
is of the form $(\ref{definition_champ})$ with a linear functional
$g$. For any $i$ in $\Z^d$, $\delta_{i,p} = \vert a_i\vert \|
\varepsilon_0 - \varepsilon^{'}_0\|_p$. So, $X$ is $p$-stable if
$$
\sum_{i\in\Z^d}\vert a_i\vert<\infty.
$$
Clearly, if $K$ is a Lipschitz continuous function, under the
above condition, the subordinated process $Y_i = K(X_i)$ is also
$p$-stable since $\delta_{i,p} = O(|a_i|)$. }
\end{example}

\begin{example}
{\em (Volterra field) Another class of nonlinear random field is
the Volterra process which plays an important role in the
nonlinear system theory (Casti \cite{Casti1985}, Rugh \cite{Rugh1981}): consider the
second order Volterra process
\begin{eqnarray*}
X_k = \sum_{s_1, s_2\in\Z^d}
 a_{s_1, s_2} \varepsilon_{k-s_1} \varepsilon_{k-s_2},
\end{eqnarray*}
where $a_{s_1, s_2}$ are real coefficients with $a_{s_1, s_2} = 0$
if $s_1 = s_2$ and $\varepsilon_i$ in $\mathbb{L}^p$, $p \geq 2$.
Let
\begin{eqnarray*}
A_k = \sum_{s_1, s_2\in\Z^d} (a_{s_1, k}^2 + a_{k, s_2}^2)
 \mbox{ and }
B_k = \sum_{s_1, s_2\in\Z^d} (|a_{s_1, k}|^p
 + |a_{k, s_2}|^p).
\end{eqnarray*}
By the Rosenthal inequality, there exists a constant $C_p > 0$
such that
\begin{eqnarray*}
\delta_{k,p} = \| X_k - X_k^*\|_p
 \le C_p A_k^{1/2} \| \varepsilon_0\|_2 \| \varepsilon_0\|_p
 + C_p B_k^{1/p} \| \varepsilon_0\|_p^2.
\end{eqnarray*}
}
\end{example}

\section{Main Results}
\label{sec:main}

To establish a central limit theorem for $S_\Gamma$ we need the
following moment inequality. With the physical dependence measure,
it turns out that the moment bound can have an elegant and concise
form.

\begin{Prop}\label{inequality}
Let $\Gamma$ be a finite subset of $\Z^d$ and $(a_i)_{i\in\Gamma}$
be a family of real numbers. For any $p\geq 2$, we have
$$
\left\|\sum_{i\in\Gamma}a_iX_i\right\|_p\leq
\left(2p\sum_{i\in\Gamma}a_i^2\right)^{\frac{1}{2}}\Delta_p
$$
where $\Delta_p=\sum_{i\in\Z^d}\delta_{i,p}$.
\end{Prop}
In the sequel, for any $i$ in $\Z^d$, we denote $\delta_i$ in
place of $\delta_{i,2}$.
\begin{Prop}\label{variance_asymptotique}
If $\Delta_2:=\sum_{i\in\Z^d}\delta_{i}<\infty$ then
$\sum_{k\in\Z^d}\vert\E(X_0X_k)\vert<\infty$. Moreover, if
$(\Gamma_n)_{n\geq 1}$ is a sequence of finite subsets of $\Z^d$
such that $\vert\Gamma_n\vert$ goes to infinity and
$|\partial\Gamma_n| / |\Gamma_n|$ goes to zero then
\begin{equation}\label{limite_variance_S_n_prime}
\lim_{n\to+\infty}\vert\Gamma_n\vert^{-1}
 \E(S^2_{\Gamma_n})=\sum_{k\in\Z^d}\E(X_0X_k).
\end{equation}
\end{Prop}

\subsection{Central Limit Theorem}

Our first main result is the following central limit theorem.
\begin{Th}\label{tlc}
Let $(X_i)_{i\in\Z^d}$ be the stationary centered random field
defined by $(\ref{definition_champ})$ satisfying $\Delta_2 :=
\sum_{i\in\Z^d}\delta_{i}<\infty$. Assume that $\sigma_n^2 :=
\E\left(S_{\Gamma_n}^2\right)\to\infty$. Let $(\Gamma_n)_{n\geq
1}$ be a sequence of finite subsets of $\Z^d$ such that
$\vert\Gamma_n\vert\to\infty$, then the Levy distance
\begin{eqnarray}
\label{eq:Levy} L[ S_{\Gamma_n}/\sqrt{|\Gamma_n|}, \,
 N(0, \sigma_n^2 / |\Gamma_n|) ] \to 0
 \mbox{ as } n \to \infty.
\end{eqnarray}
\end{Th}

We emphasize that in Theorem \ref{tlc} no condition on the domains
$\Gamma_n$ is imposed other than the natural one $|\Gamma_n| \to
\infty$. Applying Proposition \ref{variance_asymptotique}, if
$|\partial\Gamma_n| / |\Gamma_n|$ goes to zero and $\sigma^2 :=
\sum_{k\in\Z^d}\E(X_0 X_k) > 0$ then
$$
\frac{S_{\Gamma_n}}{\sqrt{\vert\Gamma_n\vert}}\converge{n}
{+\infty}{\textrm{$\mathcal{L}$}}\mathcal{N}(0,\sigma^2).
$$

Theorem \ref{tlc} can be applied to the mean estimation problem:
suppose that a stationary random field $X_i$ with unknown mean
$\mu = \E X_i$ is observed on the domain $\Gamma$. Then $\mu$ can
be estimated by the sample mean $\hat \mu = S_\Gamma / |\Gamma|$
and a confidence interval for $\hat \mu$ can be constructed if
there is a consistent estimate for ${\rm var}(S_\Gamma) /
|\Gamma|$.

Interestingly, the Theorem can also be applied to the estimation
of auto-covariance functions. For $k \in \Z^d$ let
\begin{eqnarray}
\gamma_k = {\rm cov}(X_0, X_k)
 = \E(X_0 X_k) - \mu^2.
\end{eqnarray}
Assume $X_i$ is observed over $i \in \Gamma$ and let $\Xi = \{ i
\in \Gamma: \, i+k \in \Gamma \}$. Then $\gamma_k$ can be
estimated by
\begin{eqnarray}
\hat \gamma_k = {1\over {|\Xi|}}
 \sum_{i \in \Xi}X_i X_{i+k} - \hat \mu^2.
\end{eqnarray}
To apply Theorem \ref{tlc}, we need to compute the physical
dependence measure for the process $Y_i := X_i X_{i+k}, i \in
\Z^d$. It turns out that the dependence for $Y_i$ can be easily
obtained from that of $X_i$. Note that
\begin{eqnarray*}
\delta_{i, p/2}(Y)
 &=& \| X_i X_{i+k} - X^*_i X^*_{i+k} \|_{p/2} \cr
 &\le& \| X_i X_{i+k} - X_i X^*_{i+k} \|_{p/2}
 + \| X_i X^*_{i+k} - X^*_i X^*_{i+k} \|_{p/2} \cr
 &\le& \|X_i\|_p \delta_{i+k, p}
 + \delta_{i, p} \| X^*_{i+k} \|_p
 = \|X_0\|_p(\delta_{i+k, p} + \delta_{i, p}).
\end{eqnarray*}
Hence, if $\Delta_4 = \sum_{i \in \Z^d}  \delta_{i, 4} < \infty$,
we have $\sum_{i \in \Z^d} \delta_{i, 2}(Y) < \infty$ and the
central limit theorem for $\sum_{i \in \Xi}X_i X_{i+k} / |\Xi|$
holds if $|\Xi| \to \infty$.

\subsection{Invariance Principles}

Now, we are going to see that an invariance principle holds too.
If $\A$ is a collection of Borel subsets of $[0,1]^{d}$, define
the smoothed partial sum process $\{S_{n}(A)\,;\,A\in\A\}$ by
\begin{equation}\label{process}
S_{n}(A)=\ds{\sum_{i\in\{1,...,n\}^{d}}}\,\lambda(nA\cap R_{i})X_{i}
\end{equation}
where $R_{i}=]i_{1}-1,i_{1}]\times...\times]i_{d}-1,i_{d}]$ is the unit cube with upper corner at $i$, $\lambda$ is the Lebesgue
measure on $\R^{d}$ and $X_i$ is defined by $(\ref{definition_champ})$. We equip the collection $\A$ with the
pseudo-metric $\rho$ defined for any $A,B$ in $\A$ by $\rho(A,B)=\sqrt{\lambda(A\Delta B)}$. To measure the size of $\A$
one considers the metric entropy: denote by $H(\A,\rho,\varepsilon)$ the logarithm of the smallest number
$N(\A,\rho,\varepsilon)$ of open balls of radius $\varepsilon$ with respect to $\rho$ which form a covering of $\A$. The function
$H(\A, \rho, .)$ is the entropy of the class $\A$. Let $\mathcal{C}(\A)$ be the space of continuous real functions on $\A$, equipped with the norm
$\|.\|_{\A}$ defined by $\| f\|_{\A}=\sup_{A\in\A}\vert f(A)\vert$.\\ A standard Brownian motion indexed by $\A$ is a mean zero Gaussian
process $W$ with sample paths in $\mathcal{C}(\A)$ and Cov$(W(A),W(B))=\lambda(A\cap B)$. From Dudley \cite{Dudley} we know that such a process exists if
\begin{equation}\label{entrop-metriq1}
\int_{0}^{1}\sqrt{H(\A,\rho,\varepsilon)}\,d\varepsilon<+\infty.
\end{equation}
We say that the invariance principle or functional
central limit theorem (FCLT) holds if the sequence
$\{n^{-d/2}S_{n}(A)\,;\,A\in\A\}$ converges in distribution to an
$\A$-indexed Brownian motion in the space $(\A)$. The first weak
convergence results for $\Q_{d}$-indexed partial sum processes
were established for i.i.d. random fields and for the collection
$\Q_{d}$ of lower-left quadrants in $[0,1]^{d}$, that is to say
the collection
$\{[0,t_{1}]\times\ldots\times[0,t_{d}]\,;\,(t_{1},\ldots,t_{d})\in[0,1]^{d}\}$.
They were proved by Wichura \cite{Wichu} under a finite variance
condition and earlier by Kuelbs \cite{Kuelbs} under additional
moment restrictions. When the dimension $d$ is reduced to one,
these results coincide with the original invariance principle of
Donsker \cite{Donsker}. Dedecker \cite{Dedecker2001} gave an $\mathbb{L}^{\infty}$-projective criterion for the
process $\{n^{-d/2}S_{n}(A)\,;\,A\in\A\}$ to converge in the space $\mathcal{C}(\A)$ to a mixture
of $\A$-indexed Brownian motions when the collection $\A$
satisfies only the entropy condition ($\ref{entrop-metriq1}$). This projective criterion is valid for martingale-difference
bounded random fields and provides a sufficient condition for $\phi$-mixing bounded random fields. For unbounded random fields, the result still holds 
provided that the metric entropy condition on the class $\A$ is reinforced (see \cite{Elmachkouri2002}). It is shown in \cite{Elmachkouri--Volny2003}
that the FCLT may be not valid for $p$-integrable martingale-difference random fields ($0\leq p<+\infty$) but it still holds if the conditional variances of 
the martingale-difference random field are assumed to be bounded a.s. (see \cite{Elmachkouri--Ouchti2006}). In this paper, we are going to establish 
the FCLT for random fields of the form $(\ref{definition_champ})$ (see Theorem \ref{invariance_principle}).\\
Following \cite{van-der-Vaart-Wellner}, we recall the definition of Vapnik-Chervonenkis classes ($VC$-classes) of sets: let $\mathcal{C}$ be a collection
of subsets of a set $\mathcal{X}$. An arbitrary set of $n$ points $F_n:=\{x_1,...,x_n\}$ possesses $2^n$ subsets.
Say that $\mathcal{C}$ {\em picks out} a certain subset from $F_n$ if this can be formed as a set of the form $C\cap F_n$
for a $C$ in $\mathcal{C}$. The collection $\mathcal{C}$ is said to {\em shatter} $F_n$ if each of its $2^n$ subsets can be picked out in this manner. The {\em VC-index}
$V(\mathcal{C})$ of the class $\mathcal{C}$ is the smallest $n$ for which no set of size $n$ is shattered by $\mathcal{C}$. Clearly,
the more refined $\mathcal{C}$ is, the larger is its index. Formally, we have
$$
V(\mathcal{C})=\inf\left\{n\,;\,\max_{x_1,...,x_n}\Delta_n(\mathcal{C},x_1,...,x_n)<2^n\right\}
$$
where $\Delta_n(\mathcal{C},x_1,...,x_n)=\#\left\{C\cap\{x_1,...,x_n\}\,;\,C\in\mathcal{C}\right\}$. Two classical examples of
$VC$-classes are the collection $\mathcal{Q}_d=\left\{[0,t]\,;\,t\in[0,1]^{d}\right\}$ and $\mathcal{Q}^{'}_d=\left\{[s,t]\,;\,s,t\in[0,1]^{d},\,s\leq t\right\}$
with index $d+1$ and $2d+1$ respectively (where $s\leq t$ means $s_i\leq t_i$ for any $1\leq i\leq d$).
Fore more about Vapnik-Chervonenkis classes of sets, one can refer to \cite{van-der-Vaart-Wellner}.\\
\\
Let $\beta>0$ and $h_{\beta}=\left((1-\beta)/\beta\right)^{\frac{1}{\beta}}\ind{\{0<\beta<1\}}$.
We denote by $\psi_{\beta}$ the Young function defined by $\psi_{\beta}(x)=e^{(x+h_{\beta})^{\beta}}-e^{h_{\beta}^{\beta}}$
for any $x$ in $\R^{+}$.
\begin{Th}\label{invariance_principle}
Let $(X_i)_{i\in\Z^d}$ be the stationary centered random field defined by $(\ref{definition_champ})$ and
let $\A$ be a collection of regular Borel subsets of $[0,1]^d$. Assume that one of the following condition holds:
\begin{itemize}
\item[$(i)$] The collection $\A$ is a Vapnik-Chervonenkis class with index $V$ and there exists
$p>2(V-1)$ such that $X_0$ belongs to $\mathbb{L}^p$ and $\Delta_p:=\sum_{i\in\Z^d}\delta_{i,p}<\infty$.
\item[$(ii)$] There exists $\theta>0$ and $0<q<2$ such that $\E[\exp(\theta\vert X_0\vert^{\beta(q)})]<\infty$
where $\beta(q)=2q/(2-q)$ and $\Delta_{\psi_{\beta(q)}}:=\sum_{i\in\Z^d}\delta_{i,\psi_{\beta(q)}}<\infty$ and such that the class $\A$ satisfies the condition
\begin{equation}\label{entrop-metriq2}
\int_{0}^{1}\left(H(\A,\rho,\varepsilon)\right)^{1/q}\,d\varepsilon<+\infty.
\end{equation}
\item[$(iii)$] $X_0$ belongs to $\mathbb{L^{\infty}}$, the class $\A$ satisfies the condition $(\ref{entrop-metriq1})$ and
$\Delta_{\infty}:=\sum_{i\in\Z^d}\delta_{i,\infty}<\infty$.
\end{itemize}
Then the sequence of processes $\{n^{-d/2}S_n(A)\,;\,A\in\A\}$ converges in distribution in $\mathcal{C}(\A)$
to $\sigma W$ where $W$ is a standard Brownian motion indexed by $\A$ and $\sigma^2=\sum_{k\in\Z^d}\E(X_0X_k)$.
\end{Th}
\section{Proofs}
\label{sec:proof}

{\em Proof of Proposition $\ref{inequality}$}. Let $\tau:\Z\to\Z^d$ be a bijection. For any $i\in\Z$, for any $j\in\Z^d$,
\begin{equation}\label{definition_P_i}
P_iX_j:=\E(X_j\vert\F_i)-\E(X_j\vert\F_{i-1})
\end{equation}
where $\F_i=\sigma\left(\varepsilon_{\tau(l)};l\leq i\right)$.
\begin{lemma}\label{majoration_P_X} For any $i$ in $\Z$ and any $j$ in $\Z^d$, we have $\|P_iX_j\|_p\leq\delta_{j-\tau(i),p}$.
\end{lemma}
{\em Proof of Lemma $\ref{majoration_P_X}$}.
\begin{align*}
\left\|P_iX_j\right\|_p&=\left\|\E(X_j\vert\F_i)-\E(X_j\vert\F_{i-1})\right\|_p\\
&=\left\|\E(X_0\vert T^j\F_i)-\E(X_0\vert T^j\F_{i-1})\right\|_p\quad\textrm{where $T^j\F_i=\sigma\left(\varepsilon_{\tau(l)-j};l\leq i\right)$}\\
&=\left\|\E\left(g\left((\varepsilon_{-s})_{s\in\Z^d}\right)\vert T^j\F_i\right)-
\E\left(g\left((\varepsilon_{-s})_{s\in\Z^d\backslash\{j-\tau(i)\}};\varepsilon^{'}_{\tau(i)-j}\right)\vert T^j\F_{i}\right)\right\|_p\\
&\leq \left\|g\left((\varepsilon_{-s})_{s\in\Z^d}\right)-g\left((\varepsilon_{-s})_{s\in\Z^d\backslash\{j-\tau(i)\}};\varepsilon^{'}_{\tau(i)-j}\right)\right\|_p\\
&=\left\|g\left((\varepsilon_{j-\tau(i)-s})_{s\in\Z^d}\right)-g\left((\varepsilon_{j-\tau(i)-s})_{s\in\Z^d\backslash\{j-\tau(i)\}};\varepsilon^{'}_{0}\right)\right\|_p\\
&=\left\|X_{j-\tau(i)}-X_{j-\tau(i)}^{\ast}\right\|_p\\
&=\delta_{j-\tau(i),p}.
\end{align*}
The proof of Lemma $\ref{majoration_P_X}$ is complete.\\
\\
For all $j$ in $\Z^d$,
$$
X_j=\sum_{i\in\Z}P_iX_j.
$$
Consequently,
$$
\left\|\sum_{j\in\Gamma}a_jX_j\right\|_p=\left\|\sum_{j\in \Gamma}a_j\sum_{i\in\Z}P_iX_j\right\|_p
=\left\|\sum_{i\in\Z}\sum_{j\in \Gamma}a_jP_iX_j\right\|_p.
$$
Since $\left(\sum_{j\in \Gamma}a_jP_iX_j\right)_{i\in\Z}$ is a martingale-difference sequence, by Burkholder inequality, we have
\begin{equation}
\left\|\sum_{j\in\Gamma}a_jX_j\right\|_p
\leq\left(2p\sum_{i\in\Z}\left\|\sum_{j\in \Gamma}a_jP_iX_j\right\|_p^2\right)^{\frac{1}{2}}
\leq\left(2p\sum_{i\in\Z}\left(\sum_{j\in \Gamma}\vert a_j\vert\left\|P_iX_j\right\|_p\right)^2\right)^{\frac{1}{2}}
\end{equation}
By the Cauchy-Schwarz inequality, we have
$$
\left(\sum_{j\in \Gamma}\vert a_j\vert\left\|P_iX_j\right\|_p\right)^2
\leq\left(\sum_{j\in \Gamma}a_j^2\left\|P_iX_j\right\|_p\right)\times\left(\sum_{j\in\Gamma}\|P_iX_j\|_p\right)
$$
and by Lemma $\ref{majoration_P_X}$,
$$
\sum_{j\in\Z^d}\|P_iX_j\|_p\leq\sum_{j\in\Z^d}\delta_{j-\tau(i),p}=\Delta_p.
$$
So, we obtain
$$
\left\|\sum_{j\in\Gamma}a_jX_j\right\|_p\leq\left(2p\Delta_p\sum_{j\in \Gamma}a_j^2\sum_{i\in\Z}\left\|P_iX_j\right\|_p\right)^{\frac{1}{2}}.
$$
Applying again Lemma $\ref{majoration_P_X}$, for any $j$ in $\Z^d$, we have
$$
\sum_{i\in\Z}\|P_iX_j\|_p\leq\sum_{i\in\Z}\delta_{j-\tau(i),p}=\Delta_p,
$$
Finally, we derive
$$
\left\|\sum_{j\in\Gamma}a_jX_j\right\|_p\leq\left(2p\sum_{j\in \Gamma}a_j^2\right)^{\frac{1}{2}}\Delta_p.
$$
The proof of Proposition $\ref{inequality}$ is complete.\\
\\
{\em Proof of Proposition $\ref{variance_asymptotique}$}. Let $k$ in $\Z^d$ be fixed. Since $X_k=\sum_{i\in\Z}P_iX_k$ where $P_i$ is defined by ($\ref{definition_P_i}$) and
$\E((P_iX_0)(P_jX_k))=0$ if $i\neq j$, we have
$$
\E(X_0X_k)=\sum_{i\in\Z}\E((P_iX_0)(P_iX_k)).
$$
Thus, we obtain
$$
\sum_{k\in\Z^d}\vert\E(X_0X_k)\vert\leq\sum_{i\in\Z}\|P_iX_0\|_2\sum_{k\in\Z^d}\|P_iX_k\|_2.
$$
Applying again Lemma $\ref{majoration_P_X}$, we derive $\sum_{k\in\Z^d}\vert\E(X_0X_k)\vert\leq \Delta_2^2<\infty$.\\
\\
In the other part, since $(X_k)_{k\in\Z^d}$ is stationary, we have
$$
\vert\Gamma_n\vert^{-1}\E(S^2_{\Gamma_n})=\sum_{k\in\Z^d}\vert\Gamma_n\vert^{-1}\vert \Gamma_n\cap(\Gamma_n-k)\vert\E(X_0X_k)
$$
where $\Gamma_n-k=\{i-k\,;\,i\in \Gamma_n\}$. Moreover
$$
\vert \Gamma_n\vert^{-1}\vert \Gamma_n\cap(\Gamma_n-k)\vert\vert\E(X_0X_k)\vert\leq\vert\E(X_0X_k)\vert
\quad\textrm{and}\quad
\sum_{k\in\Z^d}\vert\E(X_0X_k)\vert<\infty.
$$
Since $\lim_{n\to+\infty}\vert \Gamma_n\vert^{-1}\vert \Gamma_n\cap(\Gamma_n-k)\vert=1$, applying the Lebesgue convergence theorem, we derive
$$
\lim_{n\to+\infty}\vert \Gamma_n\vert^{-1}\E(S^2_{\Gamma_n})=\sum_{k\in\Z^d}\E(X_0X_k).
$$
The proof of Proposition $\ref{variance_asymptotique}$ is complete.\\
\\
{\em Proof of Theorem $\ref{tlc}$}. We first assume that
$\liminf_{n}{\sigma_n^2}/{\vert\Gamma_n\vert}>0$. Let $(m_n)_{n\geq 1}$ be a sequence of positive integers going to infinity. In the sequel, 
we denote $\overline{X}_j=\E\left(X_j\vert\F_{m_n}(j)\right)$ where $\F_{m_n}(j)=\sigma(\varepsilon_{j-s}\,;\,\vert s\vert\leq m_n)$.
By factorization, there exists a measurable function $h$ such that $\overline{X}_j=h(\varepsilon_{j-s}\,;\,\vert s\vert\leq m_n)$. So, we have
\begin{equation}
\overline{X}^{\ast}_j
 =h(\varepsilon^{\ast}_{j-s}\,;\,\vert s\vert\leq m_n)=\E\left(X_j^{\ast}\vert\F^{\ast}_{m_n}(j)\right)
\end{equation}
where $\F^{\ast}_{m_n}(j)=\sigma(\varepsilon^{\ast}_{j-s}\,;\,\vert s\vert\leq m_n)$. We denote also for any $j$ in $\Z^d$, 
$$
\delta^{(m_n)}_{j,p}=\left\|(X_j-\overline{X}_j)-(X_j-\overline{X}_j)^{\ast}\right\|_p.
$$
The following result is a direct consequence of Proposition $\ref{inequality}$.
\begin{Prop}\label{inequality_bis}
Let $\Gamma$ be a finite subset of $\Z^d$ and $(a_i)_{i\in\Gamma}$ be a family of real numbers. For any $n$ in $\N^{\ast}$ and any $p\in[2,+\infty]$, we have
$$
\left\|\sum_{j\in \Gamma}a_j(X_j-\overline{X}_j)\right\|_p\leq\left(2p\sum_{i\in\Gamma}a_i^2\right)^{\frac{1}{2}}\Delta_p^{(m_n)}
$$
where $\Delta_p^{(m_n)}=\sum_{j\in\Z^d}\delta^{(m_n)}_{j,p}$.
\end{Prop}
We need also the following lemma.
\begin{lemma}\label{limite_delta_j_m}
Let $p\in]0,+\infty]$ be fixed. If $\Delta_p<\infty$ then $\Delta_p^{(m_n)}\to 0$ as $n\to\infty$.
\end{lemma}
{\em Proof of Lemma $\ref{limite_delta_j_m}$}. Let $j$ in $\Z^d$ be fixed. Since $(X_j-\overline{X}_j)^{\ast}=X_j^{\ast}-\overline{X}_j^{\ast}$, we have 
\begin{align*}
\delta_{j,p}^{(m_n)}&=\left\|(X_j-\overline{X}_j)-(X_j-\overline{X}_j)^{\ast}\right\|_p\leq\|X_j-X_j^{\ast}\|_p+\|\overline{X}_j-\overline{X}_j^{\ast}\|_p\\
&=\delta_{j,p}+\|\E(X_j\vert\F_{m_n}(j)\vee\F^{\ast}_{m_n}(j))-\E(X_j^{\ast}\vert\F^{\ast}_{m_n}(j)\vee\F_{m_n}(j))\|_p\\
&\leq 2\delta_{j,p}.
\end{align*}
Moreover, $\lim_{n\to +\infty}\delta_{j,p}^{(m_n)}=0$. Finally, applying the Lebesgue convergence theorem, we obtain
$\lim_{n\to+\infty}\Delta_p^{(m_n)}=0$. The proof of Lemma $\ref{limite_delta_j_m}$ is complete.\\
\\
Let $(\Gamma_n)_{n\geq 1}$ be a sequence of finite subsets of $\Z^d$ such that $\lim_{n\to+\infty}\vert\Gamma_n\vert=\infty$ and
$\liminf_n\frac{\sigma_n^2}{\vert\Gamma_n\vert}>0$ and recall that $\Delta_2$ is assumed to be finite.
Combining Proposition $\ref{inequality_bis}$ and Lemma $\ref{limite_delta_j_m}$, we have
\begin{equation}\label{ecart_X_et_X_barre}
\limsup_{n\to+\infty}\frac{\left\|S_n-\overline{S}_n\right\|_2}{\sigma_n}=0.
\end{equation}
We are going to apply the following central limit theorem due to Heinrich (\cite{Heinrich88}, Theorem 2).
\begin{Th}[Heinrich (1988)]\label{Theoreme_Heinrich}
Let $(\Gamma_n)_{n\geq 1}$ be a sequence of finite subsets of $\Z^d$ with $\vert\Gamma_n\vert\to\infty$ as $n\to\infty$ and let $(m_n)_{n\geq 1}$ be a sequence of positive
integers. For each $n\geq 1$, let $\{U_n(j),j\in\Z^d\}$ be an $m_n$-dependent random field with
$\E U_n(j)=0$ for all $j$ in $\Z^d$. Assume that $\E\left(\sum_{j\in\Gamma_n}U_n(j)\right)^2\to\sigma^2$ as $n\to\infty$ with $\sigma^2<\infty$.
Then $\sum_{j\in\Gamma_n}U_n(j)$ converges in distribution to a Gaussian random variable
with mean zero and variance $\sigma^2$ if there exists a finite constant $c>0$ such that for any $n\geq 1$,
$$
\sum_{j\in\Gamma_n}\E U_n^2(j)\leq c
$$
and for any $\varepsilon>0$ it holds that
$$
\lim_{n\to+\infty}L_n(\varepsilon):=m_n^{2d}\sum_{j\in\Gamma_n}\E\left(U_n^2(j)\ind{\vert U_n(j)\vert\geq\varepsilon m_n^{-2d}}\right)=0.
$$
\end{Th}
Since $\liminf_{n}\frac{\sigma_n^2}{\vert\Gamma_n\vert}>0$, there exists $c_0>0$ and $n_0\in\N$ such that $\frac{\vert\Gamma_n\vert}{\sigma_n^2}\leq c_0$ for any $n\geq n_0$.
Consider $S_n=\sum_{i\in\Gamma_n}X_i$, $\overline{S}_n=\sum_{i\in\Gamma_n}\overline{X}_i$ and $U_n(j):=\frac{\overline{X}_j}{\sigma_n}$. We have 
$$
\E\left(\sum_{j\in\Gamma_n}U_n(j)\right)^2=\frac{\E(\overline{S}^2_n)-\sigma_n^2}{\sigma_n^2}+1.
$$
So, for any $n\geq n_0$ we derive
\begin{align*}
\frac{\left\vert\sigma_n^2-\E(\overline{S}^2_n)\right\vert}{\sigma_n^2}
&=\frac{1}{\sigma_n^2}\left\vert\E\left(\left(\sum_{j\in\Gamma_n}(\overline{X}_j-X_j)\right)\left(\sum_{j\in\Gamma_n}(\overline{X}_j+X_j)\right)\right)\right\vert\\
&\leq\frac{1}{\sigma_n^2}\left\|\sum_{j\in\Gamma_n}(\overline{X}_j-X_j)\right\|_2
\left\|\sum_{j\in\Gamma_n}(\overline{X}_j+X_j)\right\|_2\\
&\leq\frac{2\vert\Gamma_n\vert\Delta_2^{(m_n)}}{\sigma_n^2}\left(4\Delta_2+2\Delta_2^{(m_n)} \right)\\
&\leq4c_0\Delta_2^{(m_n)}\left(2\Delta_2+\Delta_2^{(m_n)} \right)\converge{n}{+\infty}{ }0.
\end{align*}
Consequently,
$$
\lim_{n\to+\infty}\E\left(\sum_{j\in\Gamma_n}U_n(j)\right)^2=1.
$$
Moreover, for any $n\geq n_0$,
$$
\sum_{j\in\Gamma_n}\E U_n^2(j)=\frac{\vert\Gamma_n\vert\E(\overline{X}_0^2)}{\sigma_n^2}\leq c_0\E(X_0^2)<\infty.
$$
Let $\varepsilon>0$ be fixed. We have
\begin{align*}
L_n(\varepsilon)&\leq c_0m_n^{2d}\E\left(\overline{X}_0^2\ind{\left\{\vert\overline{X}_0\vert\geq \frac{\varepsilon\sigma_n}{m_n^{2d}}\right\}}\right)
\leq c_0m_n^{2d}\E\left(X_0^2\ind{\left\{\vert\overline{X}_0\vert\geq \frac{\varepsilon\sigma_n}{m_n^{2d}}\right\}}\right)\\
&\leq c_0m_n^{2d}\sigma_n\P\left(\vert\overline{X}_0\vert\geq\frac{\varepsilon\sigma_n}{m_n^{2d}}\right)
+c_0m_n^{2d}\E\left(X_0^2\ind{\left\{\vert X_0\vert\geq\sqrt{\sigma_n}\right\}}\right)\\
&\leq\frac{c_0\E(X_0^2)m_n^{6d}}{\varepsilon^2\sigma_n}+c_0m_n^{2d}\psi(\sqrt{\sigma_n})
\end{align*}
where $\psi(x)=\E\left(X_0^2\ind{\left\{\vert X_0\vert\geq x\right\}}\right)$. 
\begin{lemma}\label{psi-psin}
If the sequence $(m_n)_{n\geq 1}$ is defined for any integer $n\geq 1$ by 
$
m_n=\min\left\{\left[\psi\left(\sqrt{\sigma_n}\right)^{\frac{-1}{4d}}\right],\left[\sigma_n^{\frac{1}{12d}}\right]\right\}
$
if $\psi(\sqrt{\sigma_n})\neq0$ and by $m_n=\left[\sigma_n^{\frac{1}{12d}}\right]$ if $\psi(\sqrt{\sigma_n})=0$ where $[\,.\,]$ is the integer part function then 
$$
m_n\to\infty,\quad\frac{m_n^{6d}}{\sigma_n}\to0\quad\textrm{and}\quad m_n^{2d}\psi\left(\sqrt{\sigma_n}\right)\to0.
$$
\end{lemma}
{\em Proof of Lemma \ref{psi-psin}}. Since $\sigma_n\to\infty$ and $\psi(\sqrt{\sigma_n})\to0$, we derive $m_n\to\infty$. Moreover,
$$
\frac{m_n^{6d}}{\sigma_n}\leq\frac{1}{\sqrt{\sigma_n}}\to0\quad\textrm{and}\quad m_n^{2d}\psi\left(\sqrt{\sigma_n}\right)\leq\sqrt{\psi\left(\sqrt{\sigma_n}\right)}\to0.
$$ 
The proof of Lemma \ref{psi-psin} is complete.\\
\\
Consequently, we obtain $\lim_{n\to\infty}L_n(\varepsilon)=0$. So, applying Theorem \ref{Theoreme_Heinrich}, we derive that
\begin{equation}\label{convergence_en_loi_Sn_barre}
\frac{\overline{S}_n}{\sigma_n}\converge{n}{+\infty}{\textrm{Law}}\mathcal{N}(0,1).
\end{equation}
Combining (\ref{ecart_X_et_X_barre}) and (\ref{convergence_en_loi_Sn_barre}), we deduce
$$
\frac{S_n}{\sigma_n}\converge{n}{+\infty}{\textrm{Law}}\mathcal{N}(0,1).
$$
Hence (\ref{eq:Levy}) holds if $\liminf_n \sigma_n^2 / |\Gamma_n| >0$. In the general case, we argue as follows: 
If (\ref{eq:Levy}) does not hold then there exists a subsequence $n'\to \infty$ such that 
\begin{equation}\label{sous_suite_n_prime}
L\left[\frac{S_{n^{'}}}{\sqrt{|\Gamma_{n^{'}}|}},\,N\left(0, \frac{\sigma^2_{n^{'}}}{|\Gamma_{n^{'}}|}\right)\right]\quad\textrm{converges to some $l$ in $]0,+\infty]$}.
\end{equation} 
Assume that $\frac{\sigma_{n^{'}}^2}{\vert\Gamma_{n^{'}}\vert}$ does not converge to zero. Then there exists a subsequence $n^{''}$ such that 
$\liminf_n\frac{\sigma_{n^{''}}^2}{\vert\Gamma_{n^{''}}\vert}>0$. By the first part of the proof of Theorem $\ref{tlc}$, we obtain 
\begin{equation}\label{sous_suite_n_seconde}
L\left[\frac{S_{n^{''}}}{\sqrt{|\Gamma_{n^{''}}|}},\,N\left(0, \frac{\sigma^2_{n^{''}}}{|\Gamma_{n^{''}}|}\right)\right]\quad\textrm{converges to $0$}.
\end{equation} 
Since (\ref{sous_suite_n_seconde}) contradicts (\ref{sous_suite_n_prime}), we have $\frac{\sigma_{n^{'}}^2}{\vert\Gamma_{n^{'}}\vert}$ converges to zero. 
Consequently $S_{n'}/\sqrt{|\Gamma_{n^{'}}|}$ converges to zero in probability and 
$L\left[\frac{S_{n^{'}}}{\sqrt{|\Gamma_{n^{'}}|}},\,N\left(0, \frac{\sigma^2_{n^{'}}}{|\Gamma_{n^{'}}|}\right)\right]$ converges to $0$ which contradicts again 
$(\ref{sous_suite_n_prime})$. Consequently, (\ref{eq:Levy}) holds. The proof of Theorem \ref{tlc} is then complete.\\

{\em Proof of Theorem $\ref{invariance_principle}$}. As usual, we have to prove the convergence of the finite-dimensional laws and the tightness of the
partial sum process $\{n^{-d/2}S_n(A)\,;\,A\in\A\}$ in $\mathcal{C}(\A)$. For any Borel subset $A$ of $[0,1]^d$, we denote by $\Gamma_{n}(A)$ the finite
subset of $\Z^{d}$ defined by $\Gamma_{n}(A)=nA\cap\Z^{d}$. We say that $A$ is a regular Borel set if $\lambda(\partial A)=0$.
\begin{Prop}\label{Proposition_type_Dedecker2001} Let $A$ be a regular Borel subset of $[0,1]^d$ with $\lambda(A)>0$. We have
$$
\lim_{n\to+\infty}\frac{\vert\Gamma_n(A)\vert}{n^d}=\lambda(A)
\quad\textrm{and}\quad
\lim_{n\to+\infty}\frac{\vert\partial\Gamma_n(A)\vert}{\vert\Gamma_n(A)\vert}=0.
$$
Moreover, if $\Delta_2$ is finite then
\begin{equation}\label{limite_approximation}
\lim_{n\to+\infty}n^{-d/2}\|S_n(A)-S_{\Gamma_n(A)}\|_2=0
\end{equation}
where $S_{\Gamma_n(A)}=\sum_{i\in\Gamma_n(A)}X_i$.
\end{Prop}
{\em Proof of Proposition $\ref{Proposition_type_Dedecker2001}$}. The first part of Proposition $\ref{Proposition_type_Dedecker2001}$ is the first part of Lemma 2 in
Dedecker \cite{Dedecker2001}. So, we are going to prove only the second part. Let $n$ be a positive integer. Arguing as in Dedecker \cite{Dedecker2001}, we have
\begin{equation}\label{decomposition}
S_n(A)-S_{\Gamma_n(A)}=\sum_{i\in W_n}a_iX_i
\end{equation}
where $a_i=\lambda(nA\cap R_i)-\ind{i\in\Gamma_n(A)}$ and $W_n$ is the set of all $i$ in $\{1,..,n\}^d$ such that $R_i\cap(nA)\neq\emptyset$ and
$R_i\cap(nA)^c\neq\emptyset$. Noting that $\vert a_i\vert\leq 1$ and applying Proposition $\ref{inequality}$, we obtain 
\begin{equation}\label{decomposition_inequality_1}
\|S_n(A)-S_{\Gamma_n(A)}\|_2
\leq2\Delta_2\,\sqrt{\sum_{i\in W_n}a_i^2}\leq2\Delta_2\sqrt{\vert W_n\vert}.
\end{equation}
Following the proof of Lemma 2 in \cite{Dedecker2001}, we have $\vert W_n\vert=o(n^d)$ and we derive $(\ref{limite_approximation})$.
The proof of Proposition $\ref{Proposition_type_Dedecker2001}$ is complete.\\
The convergence of the finite-dimensional laws follows from Proposition $\ref{Proposition_type_Dedecker2001}$ and Theorem $\ref{tlc}$.\\
\\
So, it suffices to establish the tightness property.
\begin{Prop}\label{tightness}
Assume that Assumption $(i)$, $(ii)$ or $(iii)$ in Theorem $\ref{invariance_principle}$ holds. Then for any $x>0$, we have
\begin{equation}\label{tightness_equality}
\lim_{\delta\to0}\limsup_{n\to+\infty}
\P\left(\sup_{\substack{A,B\in\A \\
\rho(A,B)<\delta}}\big\vert
n^{-d/2}S_{n}(A)-n^{-d/2}S_{n}(B)\big\vert>x\right)=0.
\end{equation}
\end{Prop}
{\em Proof of Proposition $\ref{tightness}$}. Let $A$ and $B$ be fixed in $\A$ and recall that $\rho(A,B)=\sqrt{\lambda(A\Delta B)}$. We have
$$
S_n(A)-S_n(B)=\sum_{i\in\Lambda_n}a_iX_i
$$
where $\Lambda_n=\{1,...,n\}^d$ and $a_i=\lambda(nA\cap R_i)-\lambda(nB\cap R_i)$. Applying Proposition \ref{inequality}, we have
\begin{equation}\label{Lipschitz_inequality_p}
n^{-d/2}\left\|S_n(A)-S_n(B)\right\|_p\leq\Delta_p\left(\frac{2p}{n^{d}}\sum_{i\in\Lambda_n}\lambda(n(A\Delta B)\cap R_i)\right)^{\frac12}\leq\sqrt{2p}\Delta_p\rho(A,B).
\end{equation}
Assume that Assumption $(i)$ in Theorem $\ref{invariance_principle}$ holds. Then there exists
a positive constant $K$ such that for any $0<\varepsilon<1$,
we have (see Van der Vaart and Wellner \cite{van-der-Vaart-Wellner}, Theorem 2.6.4)
$$
N(\A,\rho,\varepsilon)\leq KV(4e)^V\left(\frac{1}{\varepsilon}\right)^{2(V-1)}
$$
where $N(\A,\rho,\varepsilon)$ is the smallest number of open balls of radius $\epsilon$ with respect to $\rho$ 
which form a covering of $\A$. So, since $p>2(V-1)$, we have
\begin{equation}\label{metric_entropy_Lp}
\int_0^1\left(N(\A,\rho,\varepsilon)\right)^{\frac{1}{p}}d\varepsilon <+\infty.
\end{equation}
Combining ($\ref{Lipschitz_inequality_p}$) and ($\ref{metric_entropy_Lp}$) and
applying Theorem 11.6 in Ledoux and Talagrand \cite{Led-Tal}, we infer that the sequence $\{n^{-d/2}S_{n}(A)\,;\,A\in\A\}$ satisfies the following
property: for each positive $\epsilon$ there exists a positive real $\delta$, depending on $\epsilon$ and on the value of
the entropy integral ($\ref{metric_entropy_Lp}$) but not on $n$, such that
\begin{equation}\label{esperance_plus_petite_epsilon}
\E\left(\sup_{\stackrel{A,B\in\A}{\rho(A,B)<\delta}}\vert n^{-d/2}S_{n}(A)-n^{-d/2}S_{n}(B)\vert\right)<\epsilon.
\end{equation}
The condition ($\ref{tightness_equality}$) is then satisfied under Assumption $(i)$ in Theorem $\ref{invariance_principle}$
and the sequence of processes $\{n^{-d/2}S_{n}(A)\,;\,A\in\A\}$ is tight in $\mathcal{C}(\A)$.\\
\\
Now, we assume that Assumption $(ii)$ in Theorem $\ref{invariance_principle}$ holds. The following technical lemma can be obtained using the
expansion of the exponential function.
\begin{lemma}\label{equivalence_normes}
Let $\beta$ be a positive real number and $Z$ be a real random
variable. There exist positive universal constants $A_{\beta}$ and
$B_{\beta}$ depending only on $\beta$ such that
$$
A_{\beta}\,\sup_{p>2}\frac{\| Z\|_{p}}{p^{1/\beta}}\leq\| Z\|_{\psi_{\beta}}\leq
B_{\beta}\,\sup_{p>2}\frac{\|Z\|_{p}}{p^{1/\beta}}.
$$
\end{lemma}
Combining Lemma $\ref{equivalence_normes}$ with $(\ref{Lipschitz_inequality_p})$, for any $0<q<2$, there exists $C_q>0$ such that
\begin{equation}\label{Lipschitz_inequality_psiq}
n^{-d/2}\left\|S_n(A)-S_n(B)\right\|_{\psi_q}\leq C_q\Delta_{\psi_{\beta(q)}}\rho(A,B)
\end{equation}
where $\beta(q)=2q/(2-q)$. Applying Theorem 11.6 in Ledoux and Talagrand \cite{Led-Tal}, for each positive $\epsilon$ there exists a positive
real $\delta$, depending on $\epsilon$ and on the value of
the entropy integral ($\ref{entrop-metriq2}$) but not on $n$, such that $(\ref{esperance_plus_petite_epsilon})$ holds.
The condition ($\ref{tightness_equality}$) is then satisfied and the
process $\{n^{-d/2}S_{n}(A)\,;\,A\in\A\}$ is tight in $\mathcal{C}(\A)$.\\
\\
Finally, if Assumption $(iii)$ in Theorem \ref{invariance_principle} holds then combining Lemma $\ref{equivalence_normes}$ with $(\ref{Lipschitz_inequality_p})$,
there exists $C>0$ such that
\begin{equation}\label{Lipschitz_inequality_psi2}
\left\|n^{-d/2}S_n(A)-n^{-d/2}S_n(B)\right\|_{\psi_2}\leq C\Delta_{\infty}\rho(A,B).
\end{equation}
Applying again Theorem 11.6 in Ledoux and Talagrand \cite{Led-Tal}, we obtain the tightness of the process $\{n^{-d/2}S_{n}(A)\,;\,A\in\A\}$
in $\mathcal{C}(\A)$. The proofs of Proposition $\ref{tightness}$ and Theorem $\ref{invariance_principle}$ are complete.$\hfill\Box$\\
\\
\textbf{Acknowledgments}. The authors thank an anonymous referee for his$\backslash$her constructive comments and Olivier Durieu and Hermine Bierm\'e for pointing us a mistake 
in the first version of the proof of Theorem $\ref{tlc}$.
\bibliographystyle{plain}
\bibliography{xbib}

\end{document}